\documentclass[12pt]{article}
\usepackage{a4wide}
\usepackage{amsmath}
\usepackage{amssymb}
\usepackage{comment}
\usepackage{hyphenat}
\usepackage{multicol}

\usepackage{graphicx}
\usepackage{epsfig}
\usepackage{epstopdf}
\usepackage{color}

\newtheorem{theorem}{Theorem}[section]
\newtheorem{definition}[theorem]{Definition}

\newcommand{\qed}{\hfill$\Box$}
\newenvironment{proof}[1][]{\addvspace{.2cm} \noindent{\bf Proof: #1}}{\qed\vspace{.3cm}}

\title{Use of Enumerative Combinatorics for proving the applicability of an asymptotic stability result on discrete-time SIS epidemics in complex networks}
\author{C.\ Rodr\'iguez Lucatero \\
Departamento de Tecnolog\'{\i}as de la Informaci\'on \\ 
Universidad Aut\'onoma Metropolitana-Cuajimalpa\\
Torre III \\
Av. Vasco de Quiroga 4871 \\
Col.Santa Fe Cuajimalpa, M\'exico, D. F.\\
C.P. 05348, M\'exico \\
email:crodriguez@correo.cua.uam.mx \\
\and  L.\ Alarc\'on Ramos \\
Departamento de Matem\'aticas Aplicadas y Sistemas \\ 
Universidad Aut\'onoma Metropolitana-Cuajimalpa \\
Torre III \\
Av. Vasco de Quiroga 4871 \\
Col.Santa Fe Cuajimalpa, M\'exico, D. F. \\
C.P. 05348, M\'exico \\
email:lalarcon@correo.cua.uam.mx}
\date{}
\begin{document}
\maketitle
\begin{abstract}
In this paper, we justify by the use of Enumerative Combinatorics, that the results obtained in \cite{Alarcon1}, where is analysed the complex dynamics of an epidemic model to identify the nodes that contribute the most to the propagation process and because of that are good candidates to be controlled in the network in order to stabilize the network to reach the extinction state, is applicable in almost all the cases.
The model analysed was proposed in \cite{Gomez1} 
and results obtained in \cite{Alarcon1} implies that it is not necessary to control all nodes, but only a minimal set of nodes if the topology of the network is not regular. 
This result could be important in the spirit of considering policies of isolation or quarantine of those nodes to be controlled. Simulation results were presented in \cite{Alarcon1} for large free-scale and regular networks, that corroborate the theoretical findings. In this article we justify the applicability of the controllability result obtained in \cite{Alarcon1}  in almost all the cases by means of the use of Combinatorics.
\\
\\
{\em Mathematics Subjects Classification}: 05A16,34H20,58E25
\\
\end{abstract}


\section{Introduction} \label{Intro}
Combinatorics is the science of combinations.
It is a very important subject in the field of Discrete Mathematics.
Among many other things, it help us to conceive methods for enumerating a wide range of objects that accomplish certain property of interest in a given domain. In particular the subfield of analytic combinatorics has as goal the precise prediction of large structured combinatorial configurations under the approach of analytic methods by the use of generating functions. Analytic combinatorics is the study of finite structures built according to a finite set of rules. Generating functions is a enumerative combinatorics tool that connects discrete mathematics and continuous analysis.
One of the applications of the enumerative combinatorics is the probabilisitic method.  
His rapid development is related with the important role played by radomness in Theoretical Computer Science. The probabilistic method is one of the most simple and beautiful noncostructive proving methods that has been used in combinatorics in the last sixty years. This method was one of must important contribution of the great hungarian mathematician Paul Erd\"os \cite{Erdos1}, \cite{Alon1}. This method is used for proving the existence of some kind of mathematical object by showing that if we chose some object randomly from a specified class, the probability that the mathematical object is of a prescribed type is more than zero. The method is applied in many fields of the mathematics as can be number theory, combinatorics, graph theory, linear algebra, information theory and computer science.   

The second author of the present article, who is also the first author of \cite{Alarcon1}, in a previous result that is less general than the model of the virus spreading on complex networks presented in \cite{Alarcon1}, was interested in establishing the conditions allowing to detect those nodes of a complex network that should be controlled such that the system can be steered to a stable virus extinction state. In fact this previous model was presented as an application example of the results obtained in a section of  \cite{Alarcon1}. The criteria obtained in \cite{Alarcon1} as well as the one obtained in the mentioned previous model, can be considered as a good option when the number of nodes to be controlled with respect to the total number of nodes is small. 
The number of nodes to be controlled depends on the topology of the complex network as well as on the values of the state transition parameters in each node. If all the nodes have the same degree and the transition parameters are the same for each node, then we can not distinguish the nodes to be controlled from the ones that don't have to be controlled. In such situation we have to control all the nodes. 
As a consequence, the applicability of the criteria for the detection of the nodes to be controlled obtained in \cite{Alarcon1} is no longer a good choice given that it becomes very expensive to control the totality of the nodes in the complex network.
In \cite{Alarcon1} the authors avoid this problem giving as argument that in the literature of complex networks the most part of the topolgies that can be found in practice are not regular. Additionally the authors of \cite{Alarcon1} generalise the refered previous model assigning to each node different internal transition probabilities, as well as different rate of interaction among nodes.
The first author of the present paper thinks that the situation where each node have the same internal transition probabilities as well as the same rate of interaction among nodes can not be discarded because is the most common case observed in the behavior of the agents in
the social networks and the one in wich is based the good predictability of the models used for simulating social phenomena \cite{Chakrabarti},\cite{Galam1}, \cite{Gomez1},\cite{Axelrod1}, \cite{Gonzalez1}, \cite{Gonzalez2} and
\cite{Klemm1}.  Because of that, in the present article we want to give an argument in favor of the applicability of the result obtained in \cite{Alarcon1}, by means the use of the mathematical tool of enumerative combinatorics.
For that end we are going to assume that the transitions probabilities of each node as well the rate of interaction among nodes are the same for all the network and use an enumerative combinatoric argument to show that even in that case the criteria of detection of nodes to be controlled still being applicable.


The paper is organized as follows: In Section \ref{epidemics}, the SIS mathematical model is introduced. 
In Section \ref{Controlproblem}, we introduce the previous model that was presented as an application example in
\cite{Alarcon1}. In srction \ref{bifurcation} we present a bifurcation analysis that allow to obtain the epidemic transition threshold of the complex network.
In Section \ref{NodeSelect}, we propose a method to determine the nodes to be controlled. 
In Section \ref{applicability}, we talk about the enumerative combinatoric argument that allow us to
justify the applicability of the criteria for the selection of nodes to be controlled obtained in Section
\ref{NodeSelect} as well as the one obtained in \cite{Alarcon1}.
In section \ref{enumerating} we make an overview of the graph enumeration problem. In Subsection \ref{enumreg}
we talk about some results obtained for the enumeration of regular graphs.
In section \ref{enumerating} we talk about some results obtained for the enumeration of connected graphs and asymptotic analysis.
Finally, the main result and conclusions are presented in Section \ref{mainresult}.

\section{Epidemics spreading in Complex networks} \label{epidemics}
Many systems have been studied using structures called complex networks due to the number of its elements and their interaction. Each node in the network has associated some variables that represent the state of the node .
 The edges that connect the nodes are weighted and can be undirected or directed. The resulting graph is an undirected network or a directed network. The state in any node depends on the state and interaction intensity  
with the nodes in his neighbour. Examples of such systems are the social networks, the computer networks and the electrical systems.

Recently many researchers have concentrated their efforts in studying the problem of controllability in complex networks \cite{Liu1}, \cite{Nepusz}, \cite{Pasqualetti},  \cite{Lombardi}, \cite{Tanner}, i.e. to drive the dynamical system from any initial state to some desired state in a finite time. 

The controllability problem 
bring up many questions, and one of the most important is to determine which nodes must be controlled. 
In \cite{Liu1}
the authors propose a method to find 
the minimum set of driver nodes that control the system. This minimum set, that is equal to the number of inputs, is determined by the maximum matching in the network. 
This approach of structured control theory makes use of tools of graph theory.

In \cite{Nepusz} the authors propose to make use of the maximum matching graph theory concept, but unlike the approach proposed in \cite{Liu1}, where dynamics was considered at node level to control the system, they propose to make use of the dynamic at edges level and associate the state of the system to the state of the edges of the network.
In the latter case, the authors of \cite{Nepusz} determine the smallest set of control paths and therefore of driver nodes.


The concept of controllability of linear systems has been applied 
in the context of complex networks. In \cite{Lombardi}, \cite{Tanner}, the authors give the requirements that must be satisfied in order to achieve the Kalman's controllability rank condition \cite{Kalman}. Finally, in \cite{Pasqualetti} the authors propose a control strategy which supposes that the selection of control nodes 
must be based on the partition of the nodes of the network, 
defined in terms of the energy to control the network and propose as well a distributed control law to drive the system to a desired state. However, in all the articles mentioned in the previous paragraphs  
a node is controlled by an external input that modifies the state of the node, but in the present paper a node is controlled only if it meet  
some property. The purpose of controlling a node consist in stating conditions in order to determine which nodes are the best ones for monitoring and controlling to stabilize the system in the extinction state  for stopping the propagation of a virus. In order to verify our theoretical results we make simulations on a scale free network as well as on a regular network.
 In contrast with other authors, we consider an undirected network where each node follows the well-know Susceptible-Infected-Susceptible (SIS) model \cite{Gomez1}, so that the control mechanism brings the system to the extinction state. Of course we are interested in determining the number of nodes that will be controlled.


\section{Control Problem Statement} \label{Controlproblem}
\begin{figure}
\centering
\includegraphics[scale=0.30]{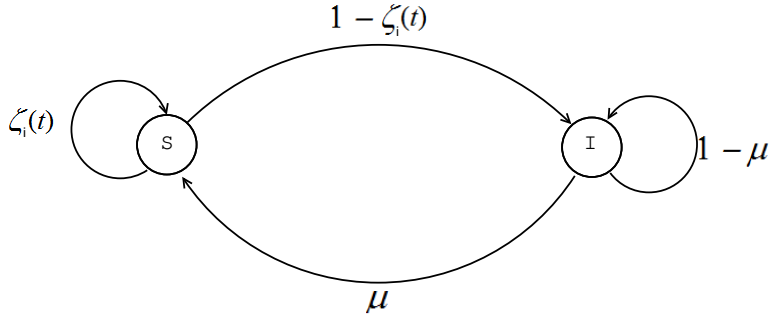}
\caption{State transition diagram associated to the SIS model, where $S$ represent the susceptible state and $I$ the infected state.}
\label{DiagEstados}
\end{figure}

The discrete time Markov process dynamical system with the SIS epidemiological model that we
will use to illustrate the problem is the one proposed in \cite{Gomez1} and was mentioned in
in a section of \cite{Alarcon1} as feedback control example.
Consider a discrete time Markov process dynamical system with the SIS epidemiological model as is described in \cite{Gomez1}.  Transitions between states depend on $\mu_{i}$ that is the recovery probability of each node $i$, and $\zeta_i(t)$ that is the probability that
a node $i$ is not being infected by 
interaction with their neighbours, as is shown in Figure \ref{DiagEstados}. In each time step the node $i$ 
try to get in   
contact with each of its neighbours, the probability that node $i$ perform at least one contact with its neighbour $j$ is given by $r_{ij}$, and the probability that the node $i$ 
will be infected by contact with an infected node is $\beta_{i}$. In that case, we consider that $\beta_{i}$ can be tuned in such a way that we can improve the nodes health. By considering  that the interaction between the node $i$ with each one of its neighbours is independent,  the probability that node $i$ is not infected by its neighbours is	

\begin{equation} \label{zetaij}
\zeta_{i}(t)=\prod_{j=1}^{N}\Big[1-a_{ij}\beta_{i} r_{i} p_{j}(t)\Big],
\end{equation} 

where $a_{ij}$ are the entries of the adjacency matrix $A \in \mathbb{R}^{N\times N}$ that represents the existing connections between the $N$ nodes of the undirected network, and $p_{i}(t)$ is the probability that a node $i$ is  infected. According to the transitions diagram (see Fig. \ref{DiagEstados}), the following non-linear dynamics is obtained:
\begin{equation}
\label{Sistema}
p_{i}(t+1)=(1-\mu_{i})p_{i}(t)+(1-\zeta_{i}(t))(1-p_{i}(t)).
\end{equation}
The purpose of the present work is to identify the nodes to be controlled in the network, so a control mechanism will allow to bring the system to the extinction state from any initial state, i.e. $p_i(t) \rightarrow  0$ as $t\rightarrow \infty$, from any state $p_i(0)=p_{oi}$, $i=1,2,\ldots N$. 
We will assume that $r_{ij}=r_{ji}=r$, 
that is, a node $i$ tries to connect to a node $ j $ in the same way that a node $ j $ tries to connect to a node $ i $.
We will assume that if a node $ i $ reduces the number of connection attempts with the node $ j $, the node $ j $ will do the same.
We have to mention that in our work we will not be considering in advance a particular complex network structure. Moreover, we will consider that not all the parameters in each node have the same value. This means that we will assume non-homogeneity both in the structure and in the properties in our model.
 
\section{Bifurcation analysis} \label{bifurcation}

Recently an important number of scientific works have been dedicated to the study of the dynamics of the spread of epidemics in complex networks
\cite{Gomez1}, \cite{Chakrabarti}, \cite{Ahn}, so we have followed the same 
research path in order to determine 
under what conditions 
the epidemic extinction state $p_i=0$ 
becomes a global attractor for any node $i$. First, we consider the following bound

\begin{equation}
\label{AproxZeta}
1-\zeta_{i}(t)\leq \sum_{j=1}^{N}a_{ij}\beta_{i} r_{i} p_{j}(t).
\end{equation}
proved in \cite{Chakrabarti} and then determine a linear dynamics substituting (\eqref{AproxZeta}) into (\eqref{Sistema})
\begin{eqnarray}
\notag
p_{i}(t+1)&=&(1-\mu_{i})p_{i}(t)+(1-p_{i}(t)) \sum_{j=1}^{N}a_{ij}\beta_{i} r_{i} p_{j}(t)\\
\notag
&\leq& (1-\mu_{i})p_{i}(t)+\sum_{j=1}^{N}a_{ij}\beta_{i} r_{i} p_{j}(t).
\end{eqnarray}
Therefore we propose the following bounded linear dynamics that is given by
\begin{equation}
\label{NuevoSistema}
x_{i}(t+1)= (1-\mu_{i})x_{i}(t)+\beta_{i} \sum_{j=1}^{N}a_{ij} r_{i} x_{j}(t), \quad p_{i}(t)\leq x_{i}(t).
\end{equation}
In matrix notation our system,   
can be summarized as $X(t+1)={\bf H}X(t)$, where ${\bf H}=[I-\bf{\mu+\beta} {\bf RA}]$, ${\bf R}=diag(r_{i})$,$\mu=diag(\mu_{i}), \beta=diag(\beta_{i})$ and $X(t)=[x_{1}(t), x_{2}(t),\ldots x_{N}(t)]^{T}$. This linear dynamics given by the matrix ${\bf H}$, is asymptotically stable if and only if the spectral radius ($\sigma(\bf H)$) of the matrix ${\bf H}$, meets the condition that $\sigma({\bf H})\in\{\lambda \in \mathbb{C}: |\lambda|<1 \}$.

If the contact probabilities $\mu_{i}$ and $\beta_{i}$  
belong to the closed interval $[\mu_{min}, \mu_{max}]$ and $[\beta_{min}, \beta_{max}]$, respectively, it is possible to bound the spectral radius as follows


\begin{equation}
\label{Threshold}
 \sigma ({\bf H})< 1.
\end{equation}

This equation is similar to that reported in \cite{Gomez1} and \cite{Chakrabarti}. Moreover, from (\eqref{Threshold}) we have several bifurcation parameters given by $\beta_{max}$, $\mu_{min}$ and $r$, and for control purposes 
it's appropriate to take $\beta$ as a control parameter. In this case, $\beta_{i}$ will be our control parameter in some nodes , in order to stabilize our system in the epidemic extinction state.

\section{Selection of nodes to be controlled} \label{NodeSelect}
For the selection of the nodes to be controlled we have to consider \eqref{NuevoSistema} as an uncoupled system due to the symmetry of the matrix  ${\bf R}$. Consider eq. (\eqref{NuevoSistema}) as a sum over all nodes as follows

The proposed model is
\begin{equation} \label{eqpi1}
p_{i}(t+1)= (1-\mu_{i})p_{i}(t)+(1-\zeta_{i}(t))(1-p_{i}(t))
\end{equation}
\begin{equation}
 \zeta_{i}(t)=\prod_{j=1}^{N}(1-a_{ij} \beta_{i} r_{i} p_{j}(t)) 
\end{equation}
It can be proved that
\begin{equation}
1-\zeta_{i}(t) \leq \sum_{j=1}^{N} a_{ij} \beta_{i} r_{i} p_{j}(t)
\end{equation}
and after substitution on \eqref{eqpi1} we obtain
\begin{equation} \label{eqpi2}
p_{i}(t+1)\leq (1-\mu_{i})p_{i}(t)+\sum_{j=1}^{N} a_{ij} \beta_{i} r_{i} p_{j}(t)
\end{equation} 
This lend us to propose the following new dynamics
\begin{equation} \label{newdynam1}
x_{i}(t+1)=(1-\mu_{i}) x_{i}(t)+\beta_{i} \sum_{j=1}^{N} a_{ij} r_{i} x_{i}(t)
\end{equation}
where $p_{i}(t) \leq x_{i}(t)$. Representing the above in matrix form we obtain the following
\begin{equation} \label{newdynam2}
\bf{X(t+1)=[ I-\mu +\beta R A] X(t)} 
\end{equation}
We can bound the expression \eqref{newdynam2} in terms of the eigenvalues, in order to determine the system stability, as follows
\begin{equation} \label{cond1}
 \sigma([ I-\mu +\beta R A] ) < 1
\end{equation}


If the condition \eqref{cond1} is satisfied the system will be asymptotically stable in the extinction state. In order to select those nodes that participate in the spreading process we have the following result known as Gerschgorin theorem \cite{Cullen1},\cite{Gerschgorin1}.
\begin{equation}\label{gersch1}
\rho(A) \leq max_{i} \sum_{j=1}^{N} |a_{ij}|
\end{equation}
which implies
\begin{equation} \label{gersch2}
|\lambda-(1-\mu_{i})| \leq \sum_{j=1}^{N} |\beta_{i} a_{ij} r_{i}|
\end{equation}
whose consequence is
\begin{equation}\label{gersch3}
|\lambda|-|1-\mu_{i}| \leq |\lambda -(1-\mu-{i})| \leq \sum_{j=1}^{N} \beta_{i} a_{ij} r_{i}
\end{equation}
from where we get
\begin{equation}\label{gersch4}
|\lambda| \leq |\lambda-1+\mu_{i}| + | 1-\mu_{i}| \leq \sum_{j=1}^{N} \beta_{i} a_{ij} r_{i} +1 -\mu_{i} \leq 1
\end{equation}
and we conclude that 
\begin{equation}\label{gersch5}
\beta_{i} \sum_{j=1}^{N}  a_{ij} r_{i} < \mu_{i}
\end{equation}

%
%


Note that the last equation does not tell us how to control the 
spread of the disease but instead tell us which nodes 
we need to control in order  to reach the epidemic extinction state. In this case, the nodes that have to be controlled are those who do not satisfy the inequality \eqref{gersch5}. 

\section{The applicability of the result} \label{applicability}
The applicability of the result obtained in the section \ref{NodeSelect} of the present article has the drawback that if the topology of the network is regular, which means that each node has the same degree, and the transition/interaction probabilities are the same for all the nodes in the network, then we cannot detect what are the nodes that have to be controlled and then we have to control all the nodes. The authors of \cite{Alarcon1} tested the applicability of a similar result by simulation under the assumption of heterogeneity of the transition probabilities
inside each node as well as in the interaction rates among the nodes. Another claim made in \cite{Alarcon1} is that 
the regular topology seldom arise in practice because the data available show that the real world computer networks 
are not homogeneous and follow a power law topology \cite{Wang1} and \cite{Barabasi}. 
The point of view of the first author of the present article is that the homogeneus case cannot be discarded because
in many mathematical models applied succesfully to social phenomena \cite{Galam1}, \cite{Galam2}, ,\cite{Axelrod1},\cite{Klemm1}
as well as to virus spreading phenomena \cite{Chakrabarti}, \cite{Chakrabarti2}, \cite{Gonzalez1}, \cite{Gonzalez2}
the model were homogeneus. The assumption of homegeneity in the aforementioned models is based on the behavior of individuals in social networks. 
Because of that, in the present article we are assuming that the model is homogeneus and we will justify under such circumstances that the result of section \ref{NodeSelect} still being applicable by proving with a combinatorial enumerative argument that the regular topology of a network is extremely unfrequent.

\section{Enumeration and generating functions} \label{enumerating}

Enumeration in Combinatorics is one of the must fertile and fascinaiting topics of Discrete Mathematics.
The purpose of this field has been to count the different ways of arranging objetcs under given constrains.
With the combinatorial enumerative methods it can be counted words, permutations, partitions, sequences and graphs.
One of the mathematical tools that is frequently used for this end are the generating functions or formal power series. The generating functions represent a bridge between discrete mathematics and continous analysis  (particulary
complex variable theory). When we face a problem whose answer is a sequence of numbers $a_{0},a_{1},a_{2}, \ldots $
and we want to know what is the closed mathematical expression that enable us to obtain the element $a_{n}$ of   that sequence sometimes we can do it at first sight by inspection. For example if the numerical sequence $1,3,5,7,9, \ldots$ we recognize at first sight that it is a sequence of odd numbers and the $n$-th element is $a_{n}= 2n-1$.
In other cases it is unreasonable to expect a simple formula as can be the case of the sequence
$2, 3, 5, 7, 11, 13, 17, 19, \ldots$, whose $a_{n}$ is the $n$-th prime number.
For some other sequences it is very hard to obtain directly a simple formula but it can be very helpful to transform it in a power series whose coefficients are the elements of that sequence as follows

\begin{equation} \label{ogf1}
\sum_{i=0}^{\infty}  a_{i} x^{i} 
\end{equation}
 
 The expression \ref{ogf1} is called {\em ordinary generating function}. This kind of series can be easly algebrically manipulated. In \cite{Segdewick2} are defined the following concepts:

\begin{definition} \label{definicion1}
A combinatorial class is a finite or denumerable set on wich a size function is defined, stisfying the folowing conditions:

\begin{enumerate}
    \item [(i)]  the size of an element is a nonnegative integer;
    \item [(ii)] the number of elements of any given size is finite.
\end{enumerate}

\end{definition}

\begin{definition} \label{definicion2}
The combinatorial classes $\mathcal{A}$ and $\mathcal{B}$ are said to be (combinatorially) isomorphic which is written
$\mathcal{A} \cong \mathcal{B}$ iff their counting sequences are identical. This condition is equivalent to the existence of a bijection from $\mathcal{A}$ to $\mathcal{B}$ that preserves size, and one also says that $\mathcal{A}$
and $\mathcal{B}$ are bijectively equivalent.
\end{definition}
 
In \cite{Segdewick2} it is mentioned that the {\em ordinary generating functions} (OGF)  as \ref{ogf1} of a sequence 
$\mathcal{A}=\{ a_{0},a_{1}, \ldots \}$ or combinarorial class $\mathcal{A}$ is the generating function of the numbers $\mathcal{A_{n}}$ whose sizes $a_{n}=card(\mathcal{A}_{n})$ such that the OGF of class $\mathcal{A}$ admits the  combinatorial form
\begin{equation} \label{combinogf1}
   A(x)= \sum_{\alpha \in \mathcal{A}} x^{|\alpha|}
\end{equation} 

It is also said that the variable $x$ marks size in the generating function. The OGF form \ref{combinogf1} can be
easly interpreted by observing that the term $x^{n}$ occurs as many times as there are objects in $\mathcal{A}$ of size $n$.

It can be defined the operation of {\bf coefficient extraction} of the term $x^{n}$ in the power series 
$A(x)=\sum  a_{n} x^{n}$ as follows:

\begin{equation} \label{coeffextract}
[x^{n}] \left(  \sum_{n \geq 0} a_{n} x^{n}  \right) = a_{n}
\end{equation}  
  
Such is the case of the sequence $0,1,1,2,3,5,8,13,21,34,55, \ldots$ know as Fibonacci numbers and that satisfy the the following recurrence relation
 
 \begin{equation} \label{Fibo1}
  F_{n+1} =  F_{n} + F_{n-1} ~~\mbox{where}~~ n \geq 1, ~~F_{0}=0, ~~F_{1}=1 
 \end{equation} 

The $n$-th element of this sequence $F_{n}$, is the coefficient of $x^{n}$ in the expansion of the function

\begin{equation} \label{Fibo2}
\frac{x}{1-x-x^{2}}
\end{equation} 

as a power series about the origin. The very interesting book \cite{Wilf1} talks about all the answers that can be
obtained by the use of the generating functions and mentioned the following list:

\begin{enumerate}

 \item Sometimes it can be found an exact formula for the members of the sequence in a pleasant way. If it is not the case, when the sequence is complicated, a good approximation can be obtained.

 \item A recurrence formula can be obtained. Most often generating functions arise from recurrence formula. Sometimes, however, a new recurrence formula, from generating functions and new insights of the nature of the sequence can be obtained.
 
 \item Averages and other statistical properties of a sequence can be obtained.
 
 \item When the sequence is very difficult to deal with asymptotic formulas can be obtained instead of an exact formula. For example for the $n$-th prime number is approximately $n \log n$ when $n$ is very big.
 
 \item Unimodality, convexity, etc. of a sequence can be proved.
 
 \item Some identities can be proved easly by using generating functions. For instance
    \begin{equation}
       \sum_{j=0}^{n} \binom{n}{j}^{2} = \binom{2n}{n} ~~(n=0,1,2,\ldots)
    \end{equation} 
   
  \item Relationship between problems can be discovered from the stricking resemblance 
        of the respective generating functions.  
     
\end{enumerate}

As was mentioned in the list of answers obtained by the use of generating functions, sometimes is hard to obtain an exact formula and in that case we can resort to the use of asymptotic formulas. The mathematical tools used for this purpose belong to field of the Analytic Combinatorics and can be consulted in \cite{Segdewick2}. The purpose of the Analytic Combinatorics is to predict with accuracy the properties of large structured combinatorial configurations with an approach based on mathematical analysis tools \cite{Segdewick2}. 
Under this approach we can start with a exact enumerative description of the combinatorial structure by the use of a generating function. This description is taken as a formal algebraic object. After that the generating function is used as an analytic object interpreting it as a mapping of the complex plane into itself. The singularities help us to determine the coefficients of the function in asymptotic form given as a result very good estimations of the counting of sequences. For this end the authors of \cite{Segdewick2} organize the Analitic Combinatorics in the following three components:

\begin{enumerate}
  \item {\em Symbolic Methods} that establish systematically relations discrete mathematics constructions and operations on generating functions that encode counting sequences.
  \item {\em Complex Asymptotics} that allow to extract asymptotic counting information from the generating functions by the mapping to the complex plane mentioned above.
  \item {\em Random structures} concerning the probabilistic properties accomplished by large random structures.
\end{enumerate}

Rich material relative to Complex Asymptotic Analysis can be found in \cite{Comtet1}. 
A nice text that can be consulted for applications of the enumerative combinatorics tools to the analysis of algorithms is \cite{Segdewick1}.

In this paper we are particularly interested in the application of the mentioned mathematical tools for the enumeration of graphs accomplishing some given properties. 
Lets define what is a graph. In the present article we will define a graph $G =<V,E>$ as a structure with a set of vertices 
$V=\{ v_{1}, v_{2},\ldots, v_{p} \}$ with cardinality $|V|=p$  called the order of $G$ 
and a set of unordered pairs of adjacent vertices, called edges,   
$E=\{\{ v_{i_{1}},v_{j_{1}} \}, \{ v_{i_{2}},v_{j_{2}} \}, \ldots, \{ v_{i_{q}},v_{j_{q}} \} \}$ if $G$ is undirected or a set of ordered pairs of adjacent vertices 
$E=\{ (v_{i_{1}},v_{j_{1}}), (v_{i_{2}},v_{j_{2}}), \ldots, (v_{i_{q}},v_{j_{q}}) \}$
if $G$ is a directed graph, having cardinality
$|E|=q$, without loops and without multiple edges. A graph $G$ with $p$ vertices and $q$ edges is called $(p,q)$ graph.
In a {\em labeled graph} of order $p$ one integer from $1$ to $p$ is assigned to each vertex as its label. 

Many questions about the number of graphs that have some specified property can be answered by the use of generating functions. Some typical questions about the number of graphs that fulfill a given property are for example: How many different graphs can I build with $n$ vertices? , How many different connected graphs with $n$ vertices exists ?, How many binary trees can be constructed with $n$ vertices ? etc.
For some of this questions the application of generating functions allow us to obtain easily a  
simple formula. For some other questions the answer is an asymptotic estimation formula.

The most commonly used generating functions are the {\em Ordinary Generating functions} and the {\em Exponential Generating functions}. 

The {\em Ordinary Generating functions} can be defined as follows
\begin{equation}
  a(x) = \sum_{k=0}^{\infty} a_{k} x^{k}
\end{equation}

where the coefficients are elements of the sequence of numbers $a_{0},a_{1},a_{2}, \ldots$
This kind of generating functions are used in combinatorial selection problems where the order is not important. 
The {\em Exponential Generating functions} can be defined as follows

\begin{equation}
  b(x) = \sum_{k=0}^{\infty} b_{k} \frac{x^{k}}{k!}
\end{equation}

where the coefficients are elements of the sequence of numbers $b_{0},b_{1},b_{2}, \ldots$
This kind of generating functions are used in combinatorial disposition problems where the order is crucial.
The sequence of numbers are counting sequences and can be encode exactly by generating functions.
Combinatorics deals with discrete objects as for example graphs, words, trees and integer partitions. One interesting problem is the enumeration of such objects. In the present article we are interested in the enumeration of graphs.
One of the main contributors of the combinatorial enumeration of graphs is the great mathematician George P\'oolya who counted graphs with given number of vertices and edges and proposed closed formulas for many graph counting problems based on group theory \cite{Polya1}. From the Polya's formulas it became relatively easy to enumerate rooted graphs, connected graphs, etc.
One of the first graph enumerating problems was the enumeration of number of triangulations of certain plane polygons
bye Leonard Euler \cite{Euler1} in the XVIII century.
After that the electrical engineer Kirchhoff in \cite{Kirchhoff1} found the number of spanning trees in connected graph, that is to say, the number of labelled trees.
Some years later the english mathematician Arthur Cayley was interested in the enumeration of trees (labeled trees, rooted trees and ordinary trees) and obtained in \cite{Cayley1} the closed formulas for solving such enumerating problem. Another great mathematician for a long time unknown and who anticipated Polya in his discoveries was J. Howard Redfield \cite{Redfield1}. Many objects and their configuration that are not graphs can be enumerated by clever transformations to graphs as for example automata, boolean functions or chemical isomers.
The generating functions are the tool used for enumerating graphs.
From the point of view of the generating functions, they are to type of graph enumerating problems:

\begin{enumerate}
  \item Labeled graphs problems
  \item Unlabeled graphs problems
\end{enumerate}

The labeled graphs problems can be easly solved with the direct application of the exponential generating functions.
The case of the unlabeled enumeration problems can be solved by using ordinary generating functions but require require the use of more combinatorial theory and the application of the P\'olya's theorem.

The first labeled graph enumeration that can be asked is: How many labeled graphs with $p$ vertices and $q$ edges are there? 
For answering that question let be $G_{p}(x)$ the polynomial or ordinary generating function whose coefficient of the term $x^{k}$ represent the number of labelled graphs with $p$ vertices and $k$ edges. If $V$ is the set of vertices of cardinality $p$, there are $q=\binom{p}{2}$ pairs of these vertices. In every vetex set $V$, each pair is adjacent or not adjacent. The number of labeled graphs with $k$ edges is therefore $\binom{q}{k}=\binom{\binom{p}{2}}{k}$. So the ordinary generating function $G_{p}(x)$ for labeled graphs with $p$ vertices is given by
\begin{equation} \label{enumlabgraph1}
   G_{p}(x) = \sum_{k=0}^{m} \binom{m}{k} x^{k}= (1+x)^{m}
\end{equation}

where $m=\binom{p}{2}$. Then the number of labeled graphs with $p$ vertices is $G_{p}(1)$ so we have that

\begin{equation} \label{enumlabgraph2}
   G_{p} = 2^{m}=2^{\binom{p}{2}}
\end{equation}

For example, if we want to know how many labeled graphs with $p=3$ vertices can be obtained we apply the formula \ref{enumlabgraph2} and we get

\begin{equation}
  G_{3}=2^{\binom{3}{2}}=2^{\frac{3!}{2! 1!}}=2^{3}=8
\end{equation}
 
If we want to know how many labeled graphs with $p=4$ vertices and exactly $q=5$ edges exist,
before the expression \ref{enumlabgraph1} we use the coefficient the term $x^{5}$

\begin{equation} \label{enumlabgraph3}
  [x^{5}] G_{p}(x) = [x^{5}]\sum_{k=0}^{6} \binom{6}{k} x^{k}=\binom{6}{5}=\frac{6!}{5! 1!}= 6 
\end{equation}

When we are working with labeled graphs it is very normal to ask in how many ways can be labeled a graph.
To  give an answer the symmetries or automorphisms of a graph have to be considered.
A graph isomorphism between a graph $G$ and a graph $G_{1}$ is a one to one map $A: V(G) \mapsto V(G_{1})$
that preserve adjacency. If $G_{1}= G$ then $A$ is called automorphism of $G$. The collection of all automorphisms of $G$, represented by $\Gamma(G)$ constitutes {\em the group} of $G$. The elements of $\Gamma(G)$ are permutations
over $V$. Let $s(G)=|\Gamma(G)|$ be the order of the group or number of symmetries of $G$. Therefore, the number
of ways in which a graph $G$ of order $p$ can be labeled is

\begin{equation} \label{groupG}
 l(G) = \frac{p!}{s(G)}
\end{equation}

Another illustrative example of graph enumeration problem is to enumerate labeled connected graphs.
A {\em path} of length $n$ can be defined as a sequence of vertices $\{ v_{0}, v_{1},\ldots,v_{n} \}$ such that the edges involved $\{v_{i},v_{i+1}\}$ for $i=0, \ldots n$ are distinct.
A {\em connected graph} is a graph in which any two vertices are joined by a {\em path}. 
If we want to obtain a formula for enumerating all the connected graphs $C_{p}$ of order $p$ we will have to resort to the notion of subgraph. 
A graph $H$ is a subgraph of a graph $G$ if $V(H) \subset V(G)$ and $E(H) \subset E(G)$. A {\em component} is a maximal connected subgraph. A {\em rooted graph} is a graph that have a distinguished vertex called {\em root}.
Two rooted graphs $H_{1}$ and $H_{2}$ are {\em isomorphic} if there exists a one to one function from $f:V(H_{1})\mapsto V(H_{2})$ that preserves the adjacency relation and the roots. A similar definition applies to the case of labelled graphs. 

Lets assume that $a_{k}$ for $k=1,2,3,\ldots$ is the number of ways of labeling  all graphs of order that fulfil the property $P(a)$ and whose exponential generating function is

\begin{equation} \label{expon1}
  a(x) = \sum_{k=1}^{\infty}  \frac{a_{k} x^{k}}{k!}
\end{equation}

Similarly, lets assume that $b_{k}$ for $k=1,2,3,\ldots$ is the number of ways of labeling  all graphs of order that fullfil the property $P(b)$ and whose exponential generating function is

\begin{equation}\label{expon2}
  b(x) = \sum_{k=1}^{\infty}  \frac{b_{k} x^{k}}{k!}
\end{equation}

If we make the product of series \ref{expon1} and \ref{expon2} the coefficients of $\frac{x^{k}}{k!}$ in $a(x)b(x)$  is the number of ordered pairs $(G_{1},G_{2})$ of two disjoint graphs where $G_{1}$ meet the property $P(a)$, $G_{2}$ fullfils the property $P(b)$, $k$ is the number of vertices in $G_{1} \cup G_{2}$ and the labels from $1$ to $k$ have been distributed over $G_{1} \cup G_{2}$.
If $C(x)$ is the exponential generating function for labeled connected graphs

\begin{equation}\label{expon3}
  C(x) = \sum_{k=1}^{\infty}  \frac{C_{k} x^{k}}{k!}
\end{equation}

then $C(x)C(x)$ is the generating function for ordered pais of labeled connected graphs. If we divide by $2$ this generating function we obtain the generating function for labeled graphs having exactly two components. If we make the same operation $n$ times we the coefficient of $\frac{x^{k}}{k!}$ represent the number of labeled graphs of order $k$ with exactly $n$ components

\begin{equation}\label{expon4}
  G(x) = \sum_{n=1}^{\infty}  \frac{C^{n}(x)}{n!}
\end{equation}

From \ref{expon4} we obtain the following relation

\begin{equation}\label{expon5}
  1 + G(x)=e^{C(x)} 
\end{equation}

Riordan in \cite{Riordan1} obtained the relation $C_{p}=J_{p}(2)$ where $J_{p}(x)$ enumerates the trees by of inversions and then deduced 

\begin{equation} \label{expon6}
  C_{p} = \sum_{k=1}^{p-1} \binom{p-2}{k-1} (2^{k}-1) C_{k} C_{p-k}
\end{equation}

From \ref{expon6} it can be noticed that if the exponential generating function for a class of graphs is known, then the exponential generating function for the class of graphs will be the logarithm of the first series, just as in \ref{expon5}.

Another equivalent recurrence that can be obtained for enumerating the connected labeled graphs of order $p$
(pag.7 in \cite{Harary1}) is:

\begin{equation} \label{expon7}
  C_{p} = 2^{\binom{p}{2}} - \frac{1}{p} \sum_{k=1}^{p-1} k \binom{p}{k} 2^{\binom{p-k}{2}} C_{k}
\end{equation}

Using \ref{expon7} the first author of the present article implemented the following Matlab code 
that can help us to obtain the number of connected graphs $C_{p}$ of orders going from $p=1$ to $p=20$
\\
\\
\begin{verbatim} 
function y = CuentaGrafConnEtiq( p )

% recurrencia que satisface el numero de grafos conexos
%  Harary Graph Enumeration pag. 7
% C_p= 2^{combinaciones(p,2)}-1/p* \sum_{k=1}^{p-1} k*
% combinaciones(p,k)*2^{combinaciones(p-k,2)}*C_{k}

C(1:p)=0;
C(1,1)=1;

for k=2:p
   acum=0;
   for j=1:k
     acum = acum + j * combinaciones(k,j) * CuentaGrafEtiq(k-j) * C(1,j);
   end
   C(1,k) = CuentaGrafEtiq(k)-(1/k)*acum;
end
y=C(1,p);
sprintf('%u',y) %% displaying the value with unsigned integer format
end

function z=combinaciones(n,k)
   z= factorial(n)/(factorial(k)*factorial(n-k));
end

%% calling the function from the matlab prompt for the calculation of the
%% evaluation from graphs of order 1 to 20

>> for i=1:20
    R(1,i)=CuentaGrafConnEtiq( i );
end    
\end{verbatim}

We get the following table

\begin{table} 
\caption{order 1 to 10}
\label{tabla1a10}
\begin{center}
\begin{tabular}{|c|c|c|c|c|c|c|c|c|c|c|}
\hline 
p & 1 & 2 & 3 & 4 & 5 & 6 & 7 & 8 & 9 & 10 \\ 
\hline 
$C_{p}$ & 1 & 1 & 4 & 38 & 728 & 26704 & 1866256 & 251548592 & 66296291072 & 34496488594816 \\ 
\hline 
\end{tabular}
\end{center}
\end{table}

\begin{table}
\caption{order 11 to 14}
\label{tabla11a14}
\begin{center}
\begin{tabular}{|c|c|c|c|c|}
\hline 
p & 11 & 12 & 13 & 14  \\ 
\hline 
$C_{p}$ & 35641657548953344 & $7.335460 \times 10^{19}$ & $3.012722 \times 10^{23}$  & $2.471649 \times 10^{27}$  \\
\hline
\end{tabular} 
\end{center}
\end{table} 

\begin{table}
\caption{order 15 to 18}
\label{tabla15a18}
\begin{center}
\begin{tabular}{|c|c|c|c|c|}
\hline 
p & 15 & 16 & 17 & 18  \\ 
\hline 
$C_{p}$ & $4.052768 \times 10^{31}$  & $1.328579 \times 10^{36}$  & $8.708969 \times 10^{40}$ & $1.41641 \times 10^{46}$  \\ 
\hline 
\end{tabular}
\end{center}
\end{table}

\begin{table}
\caption{order 19 to 20}
\label{tabla19a20}
\begin{center}
\begin{tabular}{|c|c|c|}
\hline 
p & 19 & 20 \\ 
\hline 
$C_{p}$ & $2.992930 \times 10^{51}$ & $1.569216 \times 10^{57}$  \\ 
\hline 
\end{tabular}
\end{center}
\end{table} 

As we can see in the table above the number of posible connected graphs $C_{p}$ grows very fast with respect to the number $p$ of vertices. Concerning the expressions \ref{expon6} and \ref{expon7}  it can be mentioned that they are recurrence relations instead of a closed formula. The recurrences \ref{expon6} or \ref{expon7} can be used for performing the calculation of $C_{p}$ with a computer program. Some recurrences are very hard to be solved and others cannot be solved but for obtaining an approximate value for $p$ very big, there are methods in analytic combinatorics that help us to calculate a very good approximation  of the $p-$th coefficient.
The generating functions a very important concept in combinatorial theory whose algebraic structure allow to reflect the structure of combinatorial clases. The approach taken by the analytic combinatorics is to examine the generating functions from the point of view of the mathematical analysis by assigning not only real value values to its variables but also values in the complex plane. The assignation of complex values to the variables of the generating functions
give as a consequence that the function becomes a geometric transformation of the complex plane. This type of geometrical mapping is regular ({\em holomorphic}) near the origin of the complex plane. Far from the origin some singularities appear that correspond to the abscence of smoothness of the function and give a lot of 
information about the function coefficients and their asymptotic growth.
Sometimes elementary real analysis suffices for estimating asymptotically enumerative sequences. At the next level of difficulty the generating functions still being explicit but its form dont allow to obtain the coefficients of the
series easly. In such cases the the complex plane analysis is a good option for estimating asymptotically these coefficients.

With the purpose of illustrating the notion of singularities let us take as an example the {\em ordinary generating function} of the Catalan numbers

\begin{equation} \label{ogfCatalan1}
f(x) = \frac{1}{2}  (1 -\sqrt{1-4x})
\end{equation}   

The expresion \ref{ogfCatalan1} is a compact description of the power series of the form

\begin{equation}
   (1-y)^{1/2} = 1 - \frac{1}{2}y - \frac{1}{8}y^{2} - \ldots
\end{equation}

The coefficients of the generating function associated to \ref{ogfCatalan1} have
the following explicit form

\begin{equation} \label{coeffogfCatalan1}
  f_{n} = [x^{n}] f(x) = \frac{1}{n} \binom{2n-2}{n-1}
\end{equation}

If we use Stirling formula we can obtain the following asymptotic 
approximation to \ref{coeffogfCatalan1}

\begin{equation} \label{StirlingfCatalan1}
   f_{n} \sim \lim_{n \rightarrow \infty} \frac{4^{n}}{\sqrt{\pi n^{3}}}
\end{equation}

The approximation of the kind of \ref{StirlingfCatalan1} can be derived by use of the generating functions as analytic objects.  For this end we can substitute in the power series expansion of the generating function $f(x)$ any real or complex value $\rho_{f}$ whose modulus small enough for example $\rho_{f}=4$. The graph associated with \ref{ogfCatalan1} is smooth and differentiable in the real plane and tends to the limit $\frac{1}{2}$ as $x \rightarrow (\frac{1}{4})^{-}$ but if we calculate its derivative we obtain the following function 

\begin{equation} \label{derogfCatalan1}
f(x) = \frac{1}{1 -\sqrt{1-4x}}
\end{equation}   

and it can be noticed that the derivative \ref{derogfCatalan1} becomes infinite in $\rho_f = \frac{1}{4}$. The points where the smoothness stops are called {\em singularities}. It can also be observed that the region where function \ref{ogfCatalan1} still being continuous can be extended, lets take for instance $x=-1$ 

\begin{equation} \label{evalogfCatalan1}
  f(-1)= \frac{1}{2} (1-\sqrt{5})
\end{equation}

We can proceed in a similar way and try to evaluate \ref{ogfCatalan1} with values in the complex plane whose modulus is less tham the radius of convergence of the series defined by \ref{ogfCatalan1} and observe that the ortogonal and regular grid in that can be defined in the real plane get transformed in a grid on the complex plane that preserves the angles of the curves of the grid which correspond to the property of complex differentiability and wich also is equivalent to the property of analycity.
As regards the asymptotic behavior of the coefficients $f_{n}$ of the generating function it can be observed that
it has a general asymptotic pattern composed by an exponential growth factor $A^{n}$ and a subexponential factor 
$\theta(n)$. In the case of \ref{StirlingfCatalan1} $A=4$ and $\theta(n) \sim \frac{1}{4}(\sqrt{\pi n^{3}})^{-1}$ so we can relate the exponential growth factor with the radius of convergence of the series by $A=\frac{1}{\rho_{f}}$ that is the singularity found along the positive real axis of the complex plane that in general correspond to the pole of the generating function and the subexponential part $\theta(n) = O(n^{-\frac{3}{2}})$ arises from the singularity of the square root type. This asymptotic behavior can be summarized as follows

\begin{equation} \label{coeffAsymptotic1}
  [x^{n}] f(x) = A^{n} \theta(n)
\end{equation} 

The exponential growth part of \ref{coeffAsymptotic1} is know as {\bf first principle of coefficient asymptotics} and
the subexponential growth part as {\bf second principle of coefficient asymptotics}. More general generating functions
can be addressed with complex variable theory results as can be the {\em Cauchy residue theorem} that relates global properties of a {\em meromorphic function} (its integral along closed curves) to purely local characteristics at the residues poles. An important application of the {\em Cauchy residue theorem} concerns coefficient of analytic functions. This is stated in the following theorem \cite{Segdewick2}

\begin{theorem}
 (Cauchy's coefficient formula). Let $f(z)$ be analytic in a region containing $0$ and let $\lambda$ be a simple loop around $0$ that is positively oriented. Then the coefficient $[z^{n}]f(z)$ admits the integral representation
 \begin{equation} \label{cauchy1}
   f_{n} \equiv [z^{n}]f(z)=\frac{1}{2i\pi} \int_{\lambda} f(z) \frac{dz}{z^{n+1}}
 \end{equation}
\end{theorem} 

For more details about Analytic Combinatorics we recommend to consult \cite{Segdewick2} as well as \cite{Comtet1}.


\subsection{Enumerating regular graphs} \label{enumreg}
As we mentioned in section \ref{enumerating} some enumerating problems can be solved easly using generating tools
for obtaining a closed formula. Some other problems are more hard to deal with and for obtaining a closed mathematical
expression but we can resort in such a case to the asymptotic approximation of the coefficients of the power series.
It was also mentioned in section \ref{enumerating} that there are some graph enumerating problems where the nodes are labelled and in such a case the use of the {\em exponential generating functions} is well adapted for this kind of problems. The other case of graph enumerating problems is when we are dealing with graphs whose nodes does not have an assigned label, then we can resort in such case to the Polya's enumerating method \cite{Polya1} \cite{Harary1} and the best choice is to use {\em ordinary generating functions}. It should also be mentioned that the edges of the graphs to be enumerated can be directed or undirected.

One of the seminal articles of enumerating graphs is \cite{Wright1} where was proved a fundamental theorem in the theory of random graphs on $n$ unlabelled nodes and with a given number $q$ of edges. 

In \cite{Wright1} the authors obtained a necessary and sufficient condition for relating asymptotically the number of unlabelled graphs 
with $n$ nodes and $q$ edges with the number of labelled graphs with $n$ nodes and $q$ edges.
Let $T_{nq}$ the number of different graphs with $n$ nodes and $q$ edges, $F_{nq}$ the corresponding number of labelled graphs, $N= \frac{n(n-1)}{2}$ the possible edges and
$F_{nq}=\binom{N}{q}=\frac{N!}{q!(N-q)!}$. The result obtained in \cite{Wright1} can be stated as the following theorem

\begin{theorem} \label{resultWright1}
 The necessary and sufficient conditions that
 \begin{equation}
     T_{nq} \sim \frac{F_{nq}}{n!}
 \end{equation}
as $n \rightarrow \infty$ is that
\begin{equation}
  \min(q,N-q)/n - (\log{n}/2) \rightarrow \infty
\end{equation}

\end{theorem} 

The formal result expressed in theorem \ref{resultWright1} for unlabelled graphs is a starting point on the enumeration of regular graphs because it allows to enumerate those unlabelled graphs 
that have some number of edges. In fact the author of \cite{Wright1} proved that if a graph with $|E|=E(n)$ edges, where $n$ is number of vertices or order of such a graph, has no isolated vertices or two vertices of degree $n-1$, then the number of unlabelled graphs of order $n$ and
number of edges $|E|$ divided by the number unlabelled graphs is asymptotic to $n!$.

Another interesting article on asymptotic enumeration of labelled graphs having a given degree sequence was \cite{Bender1}. The authors of \cite{Bender1} obtained their asymptotic result for
$n \times n$ symmetric matrices subject to the following constrains:
\begin{enumerate}
  \item [(i)]  each row sum is specified and bounded
  \item [(ii)] the entries are bounded
  \item [(iii)] a specified {\em sparse} set of entries must be zero
\end{enumerate}
The authors of \cite{Bender1} mentioned that their results can be interpreted in terms of incidence matrices for labelled graphs.
The results of \cite{Bender1} can be stated as follows.
Let ${\mathcal M}(n,z)$ be the set of all $n \times n$ symmetric $(0,1)$ matrices with at most $z$ zeroes in each row, ${\bf r}$ a vector over $[d]=\{0,1,\ldots d\}$ and $G(M,{\bf r},t)$ the number of $n \times n$ symmetric matrices $(g_{ij})$ over $[t]= \{0,1,\ldots t \}$ such that
\begin{enumerate}
    \item [(i)] $g_{ij}=0$ if $m_{ij}=0$ 
    \item [(ii)] $\sum_{j} g_{ij}=r_{i}$
\end{enumerate} 
 
\begin{theorem} \label{Bender1}
  For given $d,t$ and $z$
  \begin{equation}
       G(M, {\bf r},t) \sim \frac{T(f,\delta) e^{\epsilon a - b}}{\prod r_{i}!}
  \end{equation}
Uniformly for $(M, {\bf r}) \in \bigcup_{n=1}^{\infty} (\mathcal{M}(n,z) \times [0,d]^{n})$ as
$f \rightarrow \infty$ where $f=\sum_{i} r_{i}, \epsilon = -1 \mbox{~if~} t =1 \mbox{~and~} +1 \mbox{~if~} t>1$, for 
$a = \Big( \frac{\sum_{i} \binom{r_{i}}{2}}{f} \Big)^{2} +\Big( \frac{\sum_{m_{ij}}\binom{r_{i}}{2}}{f}\Big),~b= \Big( \sum_{m_{ij}=0, i<j} r_{i} r_{j} + \sum_{i} \frac{\binom{r_{i}}{2}}{f} \Big), ~\delta = \sum_{m_{ij}=0} r_{i}$  and $T(f,\delta)$ being the number of involutions on $[1,f]$ such that no element in some specified set of size $\delta$ is fixed.
  
\end{theorem} 
 
Three years later appeared the article \cite{Bollobas1} given a different approach of \cite{Bender1} allowing to obtain a more general asymptotic formula without reference to an exact formula. The asymptotic result obtained by Bela Bollobas in \cite{Bollobas1} for eumerating labelled regular graphs is proved by a probabilistic method. This result can be stated as follows.
Let $\Delta$ and $n$ be natural numbers such that $\Delta n = 2m$ is even 
and $\Delta \leq (2 \log{n})^{\frac{1}{2}}$, where $n$ is the number of vertices and $m$ is the number of edges of the graph $G$. Then as $n \rightarrow \infty$, the number of labelled $\Delta-$regular graphs on $n$ vertices is asymptotic to
\begin{equation} \label{BollobasLabelled}
 e^{-\lambda-\lambda^{2}} \frac{(2m)!}{m! 2^{m}(\Delta!)^{m}}
\end{equation} 

where $\lambda = \frac{(\Delta - 1)}{2}$.  

The authors of \cite{Bollobas1} affirm that the asymptotic formula \ref{BollobasLabelled} holds not only for $\Delta$ contant but also for $\Delta$ growing slowly as $n \rightarrow \infty$ and
summarized this in the following theorem.

\begin{theorem} \label{BollobasLabelledTheo1}
Let $d_{1} \geq d_{2} \geq \ldots d_{n}$ be natural numbers wth $\sum_{i=1}^{n} d_{i} = 2m$ even.
Suppose $\Delta = d_{1} \leq (2 \log{n})^{\frac{1}{2}}-1$ and 
$m \geq \max{\{ \epsilon \Delta n, (1+\epsilon)n  \}}$ for some $\epsilon > 0$. Then the number
$L({\bf d})$ of labelled graphs with degree sequence ${\bf d}=(d_{i})_{1}^{n}$ satisfies
\begin{equation} \label{Bollobas1res}
 L({\bf d}) \sim e^{-\lambda-\lambda^{2}} \frac{(2m)_{m}}{\{ 2^{m} \prod_{i=1}^{n} d_{i}!\}}
\end{equation}

where $\lambda=\frac{1}{2m} \sum_{i=1}^{n} \binom{d_{i}}{2}$

\end{theorem}

In the next year the author on \cite{Bollobas1} extended this result to the case of unlabelled graphs in \cite{Bollobas2}. The result of theorem \ref{BollobasLabelledTheo1} extended for the case of unlabelled graphs can be summarized in the following theorem

\begin{theorem} \label{BollobasUnlabelledTheo2}
 If $\Delta \geq 3$ and $L_{\Delta}=e^{-\lambda-\lambda^{2}} \frac{(2m)!}{m! 2^{m}(\Delta!)^{m}}$ then
\begin{equation}
 U_{\Delta} \sim \frac{L_{\Delta}}{n!} \sim e^{-\frac{(\Delta^{2}-1)}{4}} \frac{(2m)!}{2^{m}m!} \frac{(\Delta!)^{-n}}{n!}
\end{equation}
where $m=\frac{\Delta n}{2}$
\end{theorem}

For the details of the proof of theorem \ref{BollobasLabelledTheo1} and theorem \ref{BollobasUnlabelledTheo2} see \cite{Bollobas1} and \cite{Bollobas2} respectively.


\section{Combinatorial proof of Applicability of the result on control node selection} \label{mainresult}

As was mentioned in the section \ref{Intro} the case of homogeneity in the behavior of the nodes and their interaction can not be discarded given that what has been observed in the reaction of the agents in the context of social networks is that they try to minimize the conflict. Many successful models as can be for example \cite{Galam1},\cite{Axelrod1}, \cite{Chakrabarti} and \cite{Gonzalez1} base their predicting effectiveness on the homogeneity of the behaviour of the nodes and their interaction. By other side, in the section \ref{NodeSelect} of the present paper we have obtained a criteria for selecting the nodes to be controlled, but such criteria fails if we have homogeneity in the behavior of the nodes and at the same time the topology of the network is regular. Then what we want to do here is to justify the applicability of obtained in section \ref{NodeSelect} keeping the homogeneity of the nodes and try to compare the 
number regular graphs with $n$ vertices with the total of graphs that can be constructed with $n$ vertices.
For this end, based on the results on combinatorial graph enumeration mentioned on the
theorems \ref{BollobasLabelledTheo1} and \ref{BollobasUnlabelledTheo2} we can state our main result as follows.
First of all we suppose that our graph is labelled, that $G=(V,E)$ is $r-$regular with $r \geq 3$ constant and $rn = 2m$ where $n=|V|$ correspond to the number of vertices and $m=|E|$ correspond to the number of edges.  Let $L_{r}$ the number of labelled regular graphs of degree $r$ whose asymptotic value is (\cite{Bollobas1}) 

\begin{equation} 
 L_{r} \sim e^{-\frac{r^{2}-1}{4}} \frac{(2m)!}{2^{m}m!} (r!)^{n}
\end{equation}

Let $G_{n}$ the number of all possible graphs with $n$ vertices whose value is

\begin{equation}
  G_{n}= 2^{\binom{n}{2}}
\end{equation}

\begin{theorem} \label{limitTheo}
If $r \geq 3$ and $nr = 2m$ then
\begin{equation}
   \lim_{n \rightarrow \infty} \frac{L_{r}}{G_{n}} = 0
\end{equation} 

\begin{proof}
If $nr=2m$ and $r$ is constant of value $c_{1}$ then we can deduce that $m=\frac{r}{2}n=\frac{c_{1}}{2}n$ and this implies that $m=O(n)$ so lets say that $m=c_{2}n$ then

\begin{eqnarray*}
  \lim_{n \rightarrow \infty}  \frac{L_{r}}{G_{n}} = &
    \lim_{n \rightarrow \infty} \frac{e^{-\frac{c_{1}^{2}-1}{4}} \frac{(2c_{2}n)! (c_{1}!)^{-n}}{2^{c_{2}n} (c_{2}n)!}}{2^{\binom{n}{2}}} \\
    = & \lim_{n \rightarrow \infty} \frac{\frac{e^{-\frac{c_{1}^{2}-1}{4}}}{(c_{1}!)^{n}} \frac{(2c_{2}n)! }{2^{c_{2}n} (c_{2}n)!}}{2^{\frac{n(n-1)}{2}}}                
\end{eqnarray*}

applying the approximation Stirling formula 
$n! \sim \sqrt{2 \pi n} \Big( \frac{n}{e} \Big)^{n}$
  
\begin{eqnarray*}
       (c_{1}!)^{n}= & \Big(\sqrt{2 \pi c_{1}} \Big( \frac{c_{1}}{e} \Big)^{c_{1}}\Big)^{n}\\
       ~then~ & \\  
          = & \lim_{n \rightarrow \infty} \frac{\frac{e^{-\frac{c_{1}^{2}-1}{4}}}{\Big(\sqrt{2 \pi c_{1}}\Big( \frac{c_{1}}{e} \Big)^{c_{1}}\Big)^{n}} \frac{(2c_{2}n)! }{2^{c_{2}n} (c_{2}n)!}}{2^{\frac{n(n-1)}{2}}}  \\
       ~simplifying & \\
         = & \lim_{n \rightarrow \infty} \frac{\frac{1}{(\sqrt{2 \pi c_{1}})^{n}(c_{1})^{c_{1}n} e^{\frac{c_{1}^{2}+4c_{1}n+1}{4}}} \frac{(2c_{2}n)! }{2^{c_{2}n} (c_{2}n)!}}{2^{\frac{n(n-1)}{2}}} 
\end{eqnarray*}
         
applying the approximation Stirling formula

\begin{eqnarray*}
  & (2c_{2}n)! = \sqrt{2 \pi (2c_{2}n)} \Big( \frac{(2c_{2}n)}{e}\Big)^{(2c_{2}n)} \\
  \mbox{~and~} & \\ 
  & (c_{2}n)! = \sqrt{2 \pi (c_{2}n)} \Big( \frac{(c_{2}n)}{e}\Big)^{(c_{2}n)} \\
  \mbox{~we get~} & \\
         = & \lim_{n \rightarrow \infty} \frac{\frac{1}{(\sqrt{2 \pi c_{1}})^{n}(c_{1})^{c_{1}n} e^{\frac{c_{1}^{2}+4c_{1}n+1}{4}}} \frac{\sqrt{2 \pi (2c_{2}n)} \Big( \frac{(2c_{2}n)}{e}\Big)^{(2c_{2}n)} }{2^{c_{2}n} \sqrt{2 \pi (c_{2}n)} \Big( \frac{(c_{2}n)}{e}\Big)^{(c_{2}n)}}}{2^{\frac{n(n-1)}{2}}}                              
\end{eqnarray*}   

simplifying

\begin{eqnarray}\label{limitTheo2}
       =& \lim_{n \rightarrow \infty} 
 \frac {\frac{\sqrt{2}}{(\sqrt{2 \pi c_{1}})^{n} c_{1}^{c_{1}n} e^{(c_{1}^{2}+4c_{1}n+1)/4}} \big(\frac{c_{2}n}{e}\big)^{c_{2}n}}{2^{\frac{n(n-1)}{2}}} \\
\end{eqnarray}

given that $r \geq 3$, that we assumed that $r$ is a constant $c_{1}$ and that $nr=2m$ then
we have that $c_{1}=2c_{2}$ and replacing that in \ref{limitTheo2} we can express it in terms of $c_{1}$ wich is the regular degree $r$ assumed as fixed then we get

\begin{equation}\label{limitTheo3}
 \lim_{n \rightarrow \infty} 
 \frac {\frac{\sqrt{2}}{(\sqrt{2 \pi c_{1}})^{n} c_{1}^{c_{1}n} e^{(c_{1}^{2}+4c_{1}n+1)/4}} \big(\frac{c_{1}n/2}{e}\big)^{c_{1}n/2}}{2^{\frac{n(n-1)}{2}}}
\end{equation}

and replacing $c_{1}$ by $r$ in \ref{limitTheo3} we  get

\begin{equation}\label{limitTheo4}
 \lim_{n \rightarrow \infty} 
 \frac {\frac{\sqrt{2}}{(\sqrt{2 \pi r})^{n} r^{rn} e^{(r^{2}+4rn+1)/4}} \big(\frac{rn/2}{e}\big)^{rn/2}}{2^{\frac{n(n-1)}{2}}}
\end{equation}

So if the degree $r$ is constant the  $\lim_{n \rightarrow \infty} \frac{L_{r}}{G_{n}} = 0$

\end{proof}

\end{theorem}

Now our main result can be stated as a consequence of theorem \ref{limitTheo}

\begin{theorem}
If we assume that all graphs are uniformly distributed and that the nodes have homogeneous behavior then the criteria for selecting nodes to be controlled obtained in section \ref{NodeSelect} is almost always applicable.

\begin{proof}
As a consequence of theorem \ref{limitTheo} we know that the probability that a regular graph appears tends to zero as $n \rightarrow \infty$. Then the mentioned criteria is
almost always applicable.
  
\end{proof}

\end{theorem}


\end{document}